\documentclass [a4,11pt]{amsart}

\usepackage{amsmath,amssymb,amsfonts,enumerate,amsthm, amscd,}

\newcommand{\lo}{\longrightarrow}

\newtheorem{thm}{Theorem}[section]
\newtheorem{cor}[thm]{Corollary}
\newtheorem{lem}[thm]{Lemma}

\newtheorem{exam}[thm]{Example}
\newtheorem{rem}[thm]{Remark}

\def\proof{{\parindent0pt {\bf Proof.\ }}}

\theoremstyle{definition}

\theoremstyle{remark}

\theoremstyle{Definition and Notation}

\begin{document}
\bibliographystyle{amsplain}


\title[amalgamation of rings defined by bezout-like conditions]{amalgamation of rings defined by bezout-like conditions}

\author{Mohammed Kabbour}
\address{Mohammed Kabbour\\Department of Mathematics, Faculty of Science and Technology of Fez, Box 2202, University S.M. Ben Abdellah Fez, Morocco.
$$ E-mail\ address:\ mkabbour@gmail.com$$}

\author{Najib Mahdou}
\address{Najib Mahdou\\Department of Mathematics, Faculty of Science and Technology of Fez, Box 2202, University S.M. Ben Abdellah Fez, Morocco.
 $$E-mail\  address:\ mahdou@hotmail.com$$}



\begin{abstract} Let $f:A\lo B$ be a ring homomorphism and let $J$
be an ideal of $B.$ In this paper, we investigate the transfer of
notions elementary divisor ring, Hermite ring and B\'ezout ring to
the amalgamation $A\bowtie^fJ.$ We provide necessary and
sufficient conditions for $ A\bowtie^fJ$ to be an elementary
divisor ring where $A$ and $B$ are integral domains. In this case
it is shown that $ A\bowtie^fJ$ is an Hermite ring if and only it
is a B\'ezout ring. In particular, we study the transfer of the
previous notions to the amalgamated duplication of a ring $A$
along an $A-$submodule $E$ of $Q(A)$ such that $E^2\subseteq E.$
\end{abstract}

\maketitle

 \begin{section} {Introduction}
 All rings considered in this paper are assumed to be commutative,
 and have identity element and all modules are unitary.\\
 \par A ring $R$ is called an elementary divisor ring (resp.
 Hermite ring) if for every matrix $M$ over $R$ there exist non
 singular matrices $P,Q$ such that $PMQ$ $\left(\mbox{resp. }MQ\right)$ is a diagonal
 matrix (resp. triangular matrix). It proved in \cite{GH} that
 a ring $R$ is an Hermite ring if and only if for all $a,b\in R,$
 there exist $a_1,b_1,d\in R$ such that $a=a_1d,\  b=b_1d,$ and
 $Ra _1+Rb_1=R.$ A ring is a B\'ezout ring if every finitely generated
 ideal is principal. It is clear that every elementary divisor ring is
 an Hermite ring, and that every Hermite ring is a B\'ezout ring.
 Following  Kaplansky \cite{K} a ring $R$ is said to be a
 valuation ring if for any two elements in $R,$ one divides the
 other. Kaplansky proved that any valuation ring is an elementary divisor ring.\\

 \par Let $A$ and $B$ be rings, $J$ an ideal of $B$ and let $f:A\lo
 B$ be a ring homomorphism. In \cite{DFF} the amalgamation of
 $A$ with $B$ along $J$ with respect to $f$ is the sub-ring of $A\times
 B$ defined by:$$A\bowtie ^fJ=\{(a,f(a)+j)\ ;\ a\in A,j\in J\}.$$
 This construction is a generalization of the amalgamated
 duplication of a ring along an ideal introduced and studied in
 \cite{DF1}, \cite{D} and in \cite{DF2}.
 Moreover, several classical construction such as $A+xK[x]$ and
 $A+xK[[x]]$ can be studied as particular case of the
 amalgamation.   \\
 \par We denote $Q(A)$ the total ring of quotients of $A.$
 Let $E$ be an $A-$submodule of $Q(A)$ such that $E^2\subseteq E,$
  $A+E$ is a sub-ring of $Q(A)$ and $E$ is an ideal of
 $A+E.$ The amalgamated duplication of $A$ along $E:$
 $$A\bowtie  E=\{(a,a+e)\ ;\ a\in A,e\in E\}$$
 is also a particular case of the amalgamation of $A$ with $A+E$
 along $E$ with to respect $f,$ where $f:A\hookrightarrow
 A+E$ is the inclusion map. In fact, the amalgamated
 duplication of $A$ along $E$ can be studied in the frame of
 amalgamation construction. Our aim in this paper is to give a
 characterization for $A\bowtie^fJ$ to be an elementary divisor
 ring, an Hermite ring and a B\'ezout ring.

 \end{section}

\bigskip
\begin{section}{\ Main Results}
\bigskip

The set of all $n\times n$ matrices with entries from a ring $R$
will be denoted by $\mathcal{M}_n(R).$ We will let
$\mathcal{G}L_n(R)$ denote the units in $\mathcal{M}_n(R).$ Let
$A$ and $B$ be rings, for every matrix
$M=((a_{i,j},b_{i,j}))_{1\leq i,j\leq n}\in\mathcal{M}_n(A\times
B)$ we shall use the notation $M_a=(a_{i,j})_{1\leq i,j\leq n},\
M_b=(b_{i,j})_{1\leq i,j\leq n}$ and $M=M_a\times M_b.$ Let
$M,N\in\mathcal{M}_n(A\times B),$ it is easy to see that the
product $MN$ of $M$ and $N$ is giving by
$MN=~(M_aN_a)\times(M_bN_b).$\\

The following lemma will be useful to provide us many statements
in this paper.

\bigskip
\begin{lem}\label{lem2.1}
Let $A$ and $B$ a pair of integral domains, $f:A\lo B$ a ring
homomorphism and let $J$ be a proper ideal of $B.$
\begin{enumerate}
    \item [{\rm(1)}] If $A\bowtie ^fJ$ is a B\'ezout ring then $f(A)\cap J=0.$
    \item [{\rm(2)}] If $A\bowtie ^fJ$ is a B\'ezout ring and  $f$ is not
    injective then $J=0.$
\end{enumerate}
\end{lem}
\proof (1) Suppose the statement is false i.e $f(A)\cap J\neq0,$
and choose an element $a\in A$ such that $0\neq f(a)\in J.$ Then
$(0,f(a))$ is an element of $A\bowtie^fJ.$ Since $A\bowtie^fJ$ is
a B\'ezout ring the ideal generated by $(0,f(a))$ and $(a,f(a))$
is principal. Hence, there exists $(d,f(d)+j)\in A\bowtie^fJ$ such
that
$$(a,f(a))\left(A\bowtie^fJ\right)+(0,f(a))\left(A\bowtie^fJ\right)=(d,f(d)+j)\left(A\bowtie^fJ\right).$$
So, there exist $(b,f(b)+x),(c,f(c)+y),(\alpha,f(\alpha)+s)$ and
$(\beta,f(\beta)+t)$ in $A\bowtie^fJ$ such that
$$\left\{\begin{array}{c}
           (0,f(a))=(d,f(d)+j)(b,f(b)+x) \\
           (a,f(a))=(d,f(d)+j)(c,f(c)+y)\\
           (d,f(d)+j)=(0,f(a))(\alpha,f(\alpha)+s)+(a,f(a))(\beta,f(\beta)+t).\\
         \end{array}\right.$$
It follows that $d\neq0$  since $a=cd$ and $f(a)\neq0.$ Also $b=0$
since $bd=0$ and $A$ is an integral domain. From the previous
equalities we deduce that
$$f(a)=(f(d)+j)x=(f(d)+j)(f(c)+y)\mbox{ and } f(d)+j=f(a)(f(\alpha)+f(\beta)+s+t).$$
Multiplying the above equality by $x,$ we get that
$1=x(f(\alpha)+f(\beta)+s+t)$ since $B$ is an integral domain. We
conclude that $x$ is a unit, but $x\in J$ hence $J=B$ which is
absurd. We have the desired result.
\par (2) Assume that $J\neq0$  and let $0\neq u\in J.$ Since $f$ is non injective
there exists $0\neq a\in\ker f.$ From the assumption we can write
$$(a,u)\left(A\bowtie^fJ\right)+(0,u)\left(A\bowtie^fJ\right)=(d,f(d)+j)\left(A\bowtie^fJ\right)$$
for some $(d,f(d)+j)$ in $A\bowtie^fJ.$ With similar proof as in
the statement (1), we get that $J=B.$ This completes the proof of
Lemma \ref{lem2.1}. \qed
\bigskip
\begin{lem}\label{lem2.2}
The following assertions holds:\begin{enumerate}
    \item[{\rm(1)}] Let $A$ and $B$ be rings. Then $A\times B$ is an
    elementary divisor ring if and only if so $A$ and $B.$
    \item[{\rm(2)}] Let $f:A\lo B$ be a ring homomorphism and
    let $J$ be an ideal of $B.$ If $A\bowtie^fJ$ is an elementary
    divisor ring then so is $A$ and $f(A)+J.$
\end{enumerate}
\end{lem}
\proof (1) We begin by showing that if $M\in\mathcal{M}_n(A\times
B)$ then $M$ is invertible if and only if so is $M_a$ and $M_b.$
We put $M=((a_{i,j},b_{i,j}))_{1\leq i,j\leq n}.$ The determinant
of $M$ is giving by
$$\det M=\sum_{\sigma\in \mathcal{S}_n}\varepsilon(\sigma)\prod_{i=1}^n(a_{i,\sigma(i)},b_{i,\sigma(i)})$$
where $\mathcal{S}_n$ denotes the set of all permutations on $n$
letters and $\varepsilon(\sigma)$ denotes the sign of $\sigma,$
for every $\sigma\in \mathcal{S}_n.$ Thus $\det M=(\det M_a,\det
M_b).$ We say that $M$ is invertible if and only $\det M$ is a
unit. Then we have the desired result.\\
Assume that $A\times B$ is an elementary divisor ring. Let $U\in
\mathcal{M}_n(A)$ then $U\times0$ is equivalent to a diagonal
matrix $D$ with entries from $A\times B.$ There is some $P,Q\in
\mathcal{G}L_n(A\times B)$ such that $P(U\times0)Q=D.$ It follows
that $P_aUQ_a=D_a$ and so $A$ is an elementary divisor ring. By
the same way we get that $B$ is an elementary divisor ring.\\
Conversely, assume that $A$ and $B$ are elementary divisor rings
and let $M\in\mathcal{M}_n(A\times B).$ Then there exist two
invertible matrices $P_1$ and $Q_1$ (resp., $P_2$ and $Q_2$) and a
diagonal matrix $D$ (resp., $\Delta$) with entries from $A$
(resp., $B$) such that $P_1M_aQ_1=D$ (resp., $P_2M_bQ_2=\Delta$).
It follows that $$(P_1\times P_2)M(Q_1\times
Q_2)=(P_1M_aQ_1)\times(P_2M_bQ_2)=D\times\Delta,$$ which is a
diagonal matrix. From the previous part of the proof $P_1\times
P_2,\ Q_1\times Q_2\in \mathcal{G}L_n(A\times B).$ This
completes the proof of (1).\\
\par (2) Let $U=(a_{i,j})_{1\leq i,j\leq n}\in \mathcal{M}_n(A)$ and let $M$
be the matrix defined by
$M=\left((a_{i,j},f(a_{i,j}))\right)_{1\leq i,j\leq n}$ with
entries from $A\bowtie^fJ.$ We have the equality $U=M_a.$ Since
$A\bowtie^fJ$ is an elementary divisor ring $M$ is equivalent to a
diagonal matrix. From the previous part of the proof we deduce
that there exist $P$ and $Q$ in $\mathcal{G}L_n(A)$ such that
$PUQ$ is a diagonal matrix. Therefore $A$ is an elementary divisor
ring. With similar proof as in above we get that $f(A)+J$ is an
elementary divisor ring. \qed
\bigskip
\begin{rem}
Let $f:A\lo B$ be a ring homomorphism, $J$  an ideal of $B$ and
let $M\in\mathcal{M}_n(A\bowtie^fJ).$ Then $M$ is invertible if
and only if so is $M_a$ and $M_b.$
\end{rem}
\proof It is sufficient to prove that if $M_a\in
\mathcal{G}L_n(A)$ and $M_b\in \mathcal{G}L_n(B)$ then
$M_a^{-1}\times M_b^{-1}\in \mathcal{G}L_n(A\bowtie^fJ).$ Let
$(a,f(a)+j)\in A\bowtie^fJ$ which is a unit in the ring $A\times
B.$ We put $x=-f\left(a^{-1}\right)(f(a)+j)^{-1}j.$ Since $J$ is
an ideal of $B,\ x\in J.$ It is easy to get the following equality
$$\left(a^{-1},f\left(a^{-1}\right)+x\right)(a,f(a)+j)=(1,1).$$ Thus
$(a,f(a)+j)^{-1}\in A\bowtie^fJ.$ We say that $\det M$ is an
element of $A\bowtie^fJ$ which is a unit in $A\times B,$ therefore
$(\det M)^{-1}\in A\bowtie^fJ.$ Consequently,
$M^{-1}\in\mathcal{M}_n(A\bowtie^fJ).$ \qed

\bigskip
\begin{thm}\label{th2.4}
Let $A$ and $B$ a pair of integral domains, $f:A\lo B$
a ring homomorphism and let $J$ be an ideal of $B.$\\
{\rm(1)} Assume that $f$ is injective.\begin{itemize}
    \item If $J=B$ then $A\bowtie^fJ$ is an elementary divisor ring if and
    only if so is $A$ and $B.$
    \item If $J\neq B$ then $A\bowtie^fJ$ is an elementary divisor ring if and
    only if so is $f(A)+J$ and $f(A)\cap J=0.$
\end{itemize}
{\rm(2)}Assume that $f$ is not injective. Then $A\bowtie^fJ$ is an
elementary divisor ring if and only if one of the following
conditions holds:\begin{itemize}
    \item $J=0$ and $A$ is an elementary divisor ring.
    \item $J=B$ and $(A,B)$ is a pair of elementary divisor rings.
\end{itemize}

\end{thm}

\proof
(1) Two cases will be considered. \\
{\bf{case 1:}} If $J=B$ then $A\bowtie^fJ=A\times B.$ By applying
condition (1) of Lemma\ref{lem2.2}, we get that $A\bowtie^fJ$
is an elementary divisor ring if and only if so is $A$ and $B.$\\
{\bf{case 2:}} If $J\neq B$ and $A\bowtie^fJ$ is elementary
divisor ring then $f(A)\cap J=0$ by Lemma \ref{lem2.1} since every
elementary divisor ring is a B\'ezout ring. On the other hand
$f(A)+J$ is an elementary divisor ring by
Lemma \ref{lem2.2}\\
Conversely, assume that $f(A)+J$ is an elementary divisor ring and
$f(A)\cap J=0.$ We claim that the natural projection
$p_B:A\bowtie^fJ \lo f(A)+J$ $(p_B(a,f(a)+j)=f(a)+j)$ is a ring
isomorphism. Indeed, $$f(a)+j=0\ \Rightarrow\ f(a)=j=0\
\Rightarrow\ a=0$$ The conclusion is
now straightforward.\\
\par (2) Assume that $A\bowtie^fJ$ is an elementary divisor ring. By
using Lemma \ref{lem2.1}, we get that $J=0$ or $J=B.$ In the first
case $A\bowtie^fJ\simeq A,$ then $A$ is an elementary divisor
ring. In the second case $A\bowtie^fJ=A\times B.$ Hence $A$ and
$B$ are elementary divisor rings by Lemma \ref{lem2.2}. The
converse of (2) is an immediate consequence of Lemma \ref{lem2.2}.
\qed

\bigskip
Theorem \ref{th2.4} enriches the literature with a new example of
a
non valuation elementary divisor ring.\\
Let $f:A\lo B$ be a ring homomorphism and let $J$ be an ideal of
$B.$ It is easy to see that: if $A\bowtie^fJ$ is a valuation ring
then so is $A.$

\begin{exam}\label{exam2.5}\rm
Let $A$ be an elementary divisor domain which is not a valuation
ring (for instance $A=\mathbb{Z}$), and let $K$ its field of
fractions. Let $K[[x]]$ denote the ring of formal power series
over $K$ in an indeterminate $x.$ By [\cite{H}, Example1 p.161],
$A+\left(xK[[x]]\right)$ is an elementary divisor ring. We
conclude that $A\bowtie ^i\left(xK[[x]]\right),$ where $i$ is the
inclusion map of $A$ into $K[[x]],$ is an elementary divisor ring.
On the other hand $A\bowtie^i\left(xK[[x]]\right)$ is not a
valuation ring. Thus
$\mathbb{Z}\bowtie^i\left(x\mathbb{Q}[[x]]\right)$ is an
elementary divisor ring which is not a valuation ring.
\end{exam}
\bigskip
\begin{cor}
Let $A$ be an integral domain, $K$ its quotient field and let $E$
be a nonzero $A-$submodule of $K$ such that $E^2\subseteq E.$ Then
$A\bowtie E$ is an elementary divisor ring if and only if so is
$A$ and $A\subseteq E.$
\end{cor}

\proof We first prove that: Any ring $R'$ between an elementary
divisor ring $R$ and its total ring $Q(R),$ is also an elementary
divisor ring.\\ Let $M=\left(\dfrac{a_{i,j}}{d}\right)_{1\leq
i,j\leq n}\in\mathcal{M}_n(R'),$ where $a_{i,j}\in R$ for each
$1\leq i,j\leq n$ and $d$ is a nonzero divisor element of $R.$
There is some invertible matrices $P$ and $Q$ with entries from
$R$ such that $P\left(a_{i,j}\right)_{1\leq i,j\leq n}Q$ is a
diagonal matrix. Set $P\left(a_{i,j}\right)_{1\leq i,j\leq
n}Q=diag(\lambda_1,...,\lambda_n).$ Multiplying this equality by
$\dfrac{1}{d},$ we get that
$PMQ=diag\left(\dfrac{\lambda_1}{d},...,\dfrac{\lambda_n}{d}\right).$
Since $PMQ\in\mathcal{M}_n(R')$ the result follows.\\
Now suppose that $A$ is an elementary divisor ring and $A\subseteq
E.$ We have $A\bowtie E=A\times E.$ From the previous part of the
proof and condition (1) of Lemma \ref{lem2.2}, we get that
$A\bowtie E$ is an elementary divisor ring. Conversely, assume
that $A\bowtie E$ is an elementary divisor ring. We have $A\bowtie
E=A\bowtie^i~E,$ where $i:A \hookrightarrow A+E$ is the inclusion
map. By using the condition (1) of Theorem \ref{th2.4}, we obtain
the following result:
\begin{itemize}
    \item If $E=A+E$ (i.e $A\subseteq E$) then $A$ and $A+E$ are
    elementary divisor rings.
    \item Otherwise $(A+E)\cap E=0$ and $A+E$ is elementary
    divisor ring.
\end{itemize}
From the assumption $(A+E)\cap E\neq0$ since $E\subseteq A+E.$ We
conclude that $A\subseteq E$ and $A$ is an elementary divisor
ring. \qed
\bigskip
\begin{exam}\rm
Let $A$ be an integral domain and let $I$ be a nonzero ideal of
$A.$ Then $A\bowtie I$ is an elementary divisor ring if and only
if so is $A$ and $I=A.$
\bigskip
\end{exam}

\begin{lem}\label{lem2.8}
Let $A$ and $B$ be a pair of rings. Then: \\
{\rm(1)} $A\times B$ is a B\'ezout ring if and only if so is $A$
and
$B.$\\
{\rm(2)} $A\times B$ is an Hermite ring if and only if so is $A$
and $B.$
\end{lem}

\proof
 (1) Suppose that $A$ and $B$ are B\'ezout rings and let $I$ be a finitely generated ideal
 of $A\times B.$ There is some ideal $I_1$ of $A$ and  $I_2$
 of $B$ such that $I=I_1\times I_2.$ If the subset
 $\{(a_1,b_1),... ,(a_n,b_n)\}$ of $A\times B$ generate $I$ then $I_1=Aa_1+\cdots+Aa_n.$ Thus
 $I_1$ is a principal ideal of $A.$ There exists $a\in I_1$ such that
 $I_1=Aa.$ By the same way, we get that there exists $b\in I_2$ such
 that $I_2=Bb.$ We deduce that $I=(A\times B)(a,b).$ Conversely
 assume that $A\times B$ is a B\'ezout ring. Let $J_1$ be a finitely
 generated ideal of $A$ and let $J=J_1\times 0.$ Then $J$ is also finitely
 generated ideal of $A\times B,$ we get that $J$ is a principal ideal of
 $A\times B.$ Hence so is $J_1,$
 therefore $A$ is a B\'ezout  ring. Also $B$ is a B\'ezout ring since
 $ A\times B\simeq B\times A.$\\
 \par (2) Assume that $A\times B$ is an Hermite ring. Let $a,a'\in A$
 then there exist $(a_1,b_1),(a'_1,b_1'),(d,\delta)\in A\times B$
 such that $$\left\{\begin{array}{c}
           (a,0)=(a_1,b_1)(d,\delta) \\
           (a',0)=(a_1',b_1')(d,\delta)\\
           A\times B=(a_1,b_1)(A\times B)+(a_1',b_1')(A\times B).\\
         \end{array}\right.$$

 Let $(\alpha,\beta),(\alpha ',\beta ')\in A\times B$ such that
 $(\alpha,\beta)(a_1,b_1)+(\alpha',\beta')(a'_1,b_1')=(1,1)$.
It follows that $a=a_1d,\ a'=a'_1d$ and $\alpha a_1+\beta a'_1=1.$
We conclude that $A$ and $B$ is a pair of Hermite rings since
$A\times B\simeq B\times A.$ The converse of the statement is
obvious. \qed

\bigskip
\begin{thm}\label{th2.9}
Let $A$ and $B$ be a pair of integral domains, $J$ an ideal of $B$
and let $f:A\lo B$ be an injective ring homomorphism. Then the
following properties are equivalent:\begin{enumerate}
    \item [{\rm(1)}] $A\bowtie^fJ$ is an Hermite ring.
    \item [{\rm(2)}] $A\bowtie^fJ$ is a B\'ezout ring.
    \item [{\rm(3)}] One of the following conditions holds:
    \begin{itemize}
        \item $J=B,\ A$ and $B$ are B\'ezout rings.
        \item $J\neq B,\ f(A)\cap J=0$ and $f(A)+J$ is a B\'ezout
        ring.
    \end{itemize}
\end{enumerate}
\end{thm}
\proof (1) $\Rightarrow$ (2): Clear.\\
(2) $\Rightarrow$ (3): Assume that $J\neq B.$ By Lemma
\ref{lem2.1} $f(A)\cap J=0.$ Then the natural projection
$p_B:A\bowtie^fJ\lo f(A)+J;\ (p_B(a,f(a)+j)=f(a)+j)$ is a ring
isomorphism since $f$ is injective. Therefore $f(A)+J$ is a
B\'ezout ring. If $J=B$ then $A$ and $B$ are B\'ezout rings by the
condition (1) of Lemma \ref{lem2.8} since
$A\bowtie^fJ=A\times B.$\\
(3) $\Rightarrow (1)$: If $J=B$ then $A$ and $B$ are Hermite rings
since every B\'ezout domain is an Hermite ring. Hence
$A\bowtie^fJ=A\times B$ is an Hermite ring. Now we assume that
$J\neq B.$ Then $A\bowtie^fJ\simeq f(A)+J$ and so $A\bowtie^fJ$ is
a B\'ezout domain. This completes the proof of Theorem
\ref{th2.9}.\qed

\bigskip

\begin{exam}\label{exam2.10}\rm
Let $A$ be a B\'ezout domain, $K$ its quotient field, and let
$K[[x]]$ denote the ring of formal power series over $K$ in an
indeterminate $x.$ Then $A\bowtie^i\left(xK[[x]]\right),$ where
$i:A\hookrightarrow K[[x]]$ is the inclusion map, is an Hermite
ring.
\end{exam}

\proof Let $\displaystyle f=\sum_{n=0}^{\infty}a_nx^n,\
g=\sum_{n=0}^{\infty}b_nx^n$ be nonzero elements of
$R=A+\left(xK[[x]]\right),$ and let $p$ (resp. $q$) denote the
least integer such that $a_p\neq0$ (resp., $b_q\neq0$). We can
write $f=a_px^p(1+xf_1)$ and $g=b_qx^q(1+xg_1),$ where $f_1,g_1\in
K[[x]].$ Since $1+xf_1,\ 1+xg_1$ are units of $R,$$$
fR+gR=a_px^pR+b_qx^qR.$$ If $p<q$ (resp., $q<p$) then
$fR+gR=a_px^pR$ (resp., $b_qx^qR$). Suppose that $p=q$ and write
$a_p=\dfrac{a}{d}$ and $b_q=\dfrac{b}{d}$ for some nonzero
elements $a,b,d$ of $A$ (where $d=1$ if $p=q=0$). By the
assumption there exist $c,a',b'\in A$ such that $a=a'c,\ b=b'c$
and $a'A+b'A=A.$ It is easy to get that $fR+gR=\dfrac{c}{d}x^pR.$
This completes the proof that $A\bowtie^i\left(xK[[x]]\right)$ is
an Hermite ring.\qed
\bigskip

\begin{cor}
Let $A$ be an integral domain, $K$ its quotient field and let $E$
be a nonzero $A-$submodule of $K$ such that $E^2\subseteq E.$ Then
the following statements are equivalent:\begin{enumerate}
    \item[{\rm(1)}] $A\bowtie E$ is an Hermite ring.
    \item [{\rm(2)}] $A\bowtie E$ is a B\'ezout ring.
    \item [{\rm(3)}] $A$ is a B\'ezout ring and $A\subseteq E.$
\end{enumerate}
\end{cor}
\proof (2) $\Rightarrow$ (3): Let $0\neq \dfrac{a}{b}\in E$. Then
$0\neq a\in A\cap E$ and so $A\cap E\neq0.$ By applying Theorem
\ref{th2.9}, we get that $A+E=E$ and  $(A,A+E)$ is a pair of
B\'ezout rings.
It follows that $A\subseteq E$ and $A$ is a B\'ezout ring.\\
(3) $\Rightarrow$ (1): By applying Lemma \ref{lem2.8} and the
condition (3) of Theorem \ref{th2.9} it is sufficient to prove
that every ring between a B\'ezout domain and its quotient field
is also B\'ezout domain. Let $R$ be a B\'ezout domain and let $R'$
be a ring such that$R\subseteq R'\subseteq qf(R).$ Let
$\dfrac{a}{d},\dfrac{b}{d}\in R'$ then we can write
$$\left\{\begin{array}{c}
           a=a'c \\
           b=b'c\\
           \alpha a'+\beta b'=1\\
         \end{array}\right.$$
for some elements $a',b',c,\alpha,\beta$ in $R.$ Hence
$\dfrac{c}{d}=\alpha \dfrac{a}{d}+\beta\dfrac{b}{d}$ is an element
of $R'.$ Thus $\dfrac{c}{d}\in R'$ and
$R'\dfrac{a}{d}+R'\dfrac{b}{d}\subseteq R\dfrac{c}{d}.$ On the
other hand, we have:  $$\frac{c}{d}\in
R\frac{a}{d}+R\frac{b}{d}\subseteq R'\frac{a}{d}+R'\frac{b}{d}.$$
It follows that $R'\dfrac{a}{d}+R'\dfrac{b}{d}= R'\dfrac{c}{d}.$
Finally, $R'$ is a B\'ezout domain. \qed

\bigskip

\begin{exam}\rm
Let $A$ be an integral domain and let $I$ be a nonzero ideal of
$A.$ Then $A\bowtie I$ is a B\'ezout ring if and only if so is $A$
and $I=A.$
\end{exam}

\bigskip

\begin{thm}
Let $A$ and $B$ be a pair of integral domains, $J$ an ideal of $B$
and let $f:A\lo B$ be a non injective ring homomorphism. Then the
following statements are equivalent:\begin{enumerate}
    \item [{\rm(1)}] $A\bowtie^fJ$ is an Hermite ring.
    \item [{\rm(2)}] $A\bowtie^fJ$ is a B\'ezout ring.
    \item [{\rm(3)}] One of the following conditions holds:\begin{itemize}
        \item $J=B,\ A$ and $B$ are B\'ezout rings.
        \item $J=0,$ and $A$ is a B\'ezout ring.
    \end{itemize}
\end{enumerate}
\end{thm}
\proof (2) $\Rightarrow$ (3): By applying condition (2) of Lemma
\ref{lem2.1}, we get that $J=0$ or $J=B.$ If $J=0$ then
$A\bowtie^fJ\simeq A$ otherwise $A\bowtie^fJ=A\times B.$ By using
Lemma \ref{lem2.8}, we have the desired implication.\\
(3) $\Rightarrow$ (1): If $J=0$ then $A\bowtie^fJ\simeq A$ and $A$
is an Hermite ring (since $A$ is an integral domain). If $J=B$
then $A\bowtie^fJ=A\times B$ is an Hermite ring by condition (2)
of Lemma \ref{lem2.8}.\qed

\end{section}

\bigskip




\bigskip\bigskip


\end{document}